\documentstyle{amsart}

\advance \topmargin by -\headheight
\advance \topmargin by -\headsep

\evensidemargin \oddsidemargin
\newtheorem{thm}{Theorem}[section]
\newtheorem{cor}[thm]{Corollary}
\newtheorem{lem}[thm]{Lemma}
\newtheorem{prop}[thm]{Proposition}
\newtheorem{conj}[thm]{Conjecture}

\numberwithin{equation}{section}

\newcommand{\thmref}[1]{Theorem~\ref{#1}}
\newcommand{\propref}[1]{Proposition~\ref{#1}}
\newcommand{\corref}[1]{Corollary~\ref{#1}}
\newcommand{\secref}[1]{\S\ref{#1}}
\newcommand{\lemref}[1]{Lemma~\ref{#1}}

\def\Q  {{\Bbb Q}}
\def\Z  {{\Bbb Z}}

\def\wfor {{w_n\hspace{-.42 cm}^{^{_{_{^{_{\rightarrow}}}}}}\hspace{0.09cm}}}
\def\wbac {{w_n\hspace{-.43 cm}^{^{_{_{^{_{\leftarrow}}}}}}\hspace{0.09cm}}}
\def\afor {{a_n\hspace{-.37 cm}^{^{_{_{^{_{\rightarrow}}}}}}\hspace{0.03cm}}}
\def\abac {{a_n\hspace{-.38 cm}^{^{_{_{^{_{\leftarrow}}}}}}\hspace{0.03cm}}}
\def\sfor {{s_n\hspace{-.37 cm}^{^{_{_{^{_{\rightarrow}}}}}}\hspace{0.03cm}}}

\def\safor {{s_{n-1}\hspace{-.72cm}^{^{_{_{^{_{\rightarrow}}}}}}\hspace{0.36cm}}}
\def\sbfor {{s_{n-2}\hspace{-.73cm}^{^{_{_{^{_{\rightarrow}}}}}}\hspace{0.38cm}}}
\def\sbbac {{s_{n-2}\hspace{-.73cm}^{^{_{_{^{_{\leftarrow}}}}}}\hspace{0.38cm}}}
\def\scfor {{s_{n-3}\hspace{-.73cm}^{^{_{_{^{_{\rightarrow}}}}}}\hspace{0.38cm}}}
\def\scbac {{s_{n-3}\hspace{-.73cm}^{^{_{_{^{_{\leftarrow}}}}}}\hspace{0.38cm}}}
\def\sdfor {{s_{n-4}\hspace{-.73cm}^{^{_{_{^{_{\rightarrow}}}}}}\hspace{0.38cm}}}
\def\sdbac {{s_{n-4}\hspace{-.73cm}^{^{_{_{^{_{\leftarrow}}}}}}\hspace{0.38cm}}}
\def\smfor {{s_{-1}\hspace{-.55cm}^{^{_{_{^{_{\rightarrow}}}}}}\hspace{0.19cm}}}
\def\smbac {{s_{-1}\hspace{-.56cm}^{^{_{_{^{_{\leftarrow}}}}}}\hspace{0.19cm}}}

\begin{document}

\title[Symmetry and Specializability]{Symmetry and Specializability in
Continued Fractions} 
\author{Henry Cohn}
\address{Department of Mathematics, Harvard University}
\email{cohn@@math.harvard.edu}
\date{October 18, 1995}

\maketitle

\section{Introduction}

A number defined by a series does not in general have an interesting
continued fraction expansion.  There are, however, some exceptions,
such as the series
$$
\sum_{n=0}^\infty\frac{1}{2^{2^n}} =
[0,1,4,2,4,4,6,4,2,4,6,2,4,6,4,4,2,\dots], 
$$
dealt with (independently) in \cite{kmosek} and \cite{shallit1}.  The
continued fraction expansion for this series has a type of symmetry
known as folding symmetry.  In this paper, we generalize folding
symmetry, and give examples such as
\begin{eqnarray*}
\sum_{n=0}^\infty\frac{1}{T_{4^n}(2)} &=& [0, 1, 1, 23, 1, 2, 1,
18815, 3, 1, 23, 3, 1, 23, 1, 2, 1, 106597754640383, \\
&& 3, 1, 23, 1, 3, 23, 1, 3, 18815, 1, 2, 1, 23, 3, 1, 23, 1, 2, 1,
18815, 3, 1, 23, \\
\phantom{\sum_{n=0}^\infty\frac{1}{T_{4^n}(2)}} && 3, 1, 23, 1, 2, 1, \ldots ],
\end{eqnarray*}
where $T_\ell(x)$ is the $\ell$-th Chebyshev polynomial.

In addition, we prove a general characterization of this sort of
series.  These series are of the form
$$\sum_{n=0}^\infty\frac{1}{f^n(m)},$$ where $m \in \Z$, $f(x) \in
\Z[x]$, and $f^n(x)$ denotes the $n$-th iterate of $f(x)$.  (For
example, in the first case $f(x) = x^2$, and in the second $f(x) =
T_4(x)$.)  In addition, they are specialized, in the terminology of
\cite{special}.  For example, the first is a special case of
$$\sum_{n=0}^\infty\frac{1}{x^{2^n}} =
[0,x-1,x+2,x,x,x-2,x,x+2,x,x-2,x+2,\dots],$$ with $x = 2$.  A
continued fraction over $\Q(x)$, such as this one, with the property
that each partial quotient has integer coefficients, is called
specializable, because when one specializes by choosing an integer
value for $x$, one gets immediately a continued fraction whose partial
quotients are integers.  The continued fraction one obtains is then
called specialized.  We prove a theorem (\thmref{done})
that determines all sums of
the form $$\sum_{n=0}^\infty\frac{1}{f^n(x)}$$ (for $f(x) \in \Z[x]$)
that have specializable continued fractions.

Before proceeding any further, we quickly review the basics of
continued fractions, and give a short account of folding symmetry.
Consider the continued fraction
$$[a_0,a_1,\dots,a_n] = a_0+
  {1\over\displaystyle a_1+
   {\strut 1\over\displaystyle a_2+
{\strut 1 \over\displaystyle {\cdots+
    {\strut 1 \over\displaystyle a_{n-1} +
     {\strut 1\over a_n}}}}}}.$$
We define $p_0 = a_0$, $q_0=1$, $p_1=a_0a_1+1$, $q_1=a_1$, and for 
$n \ge 2$,
$$p_n = a_np_{n-1}+p_{n-2} \qquad \hbox{and} \qquad q_n =
a_nq_{n-1}+q_{n-2}.$$

\begin{prop}
\label{basics}
For each $n$,
$$
\frac{p_n}{q_n} = [a_0,a_1,\dots,a_n],$$
and we have $p_nq_{n-1} - q_np_{n-1} = (-1)^{n-1}.$
\end{prop}

For a proof of this standard result, see Theorems 149 and 150
of \cite{hardy-wright}.
Now consider the function $f_n(z)$ of $z$ defined by
\begin{equation*}
\label{ratdef}
f_n(z) = a_0+
  {1\over\displaystyle a_1+
   {\strut 1\over\displaystyle a_2+
{\strut 1 \over\displaystyle {\cdots+
    {\strut 1 \over\displaystyle a_n +
     {\strut 1\over z}}}}}}.
\end{equation*}
This function will be useful shortly.  By \propref{basics},
$$
f_n(z) = \frac{p_nz+p_{n-1}}{q_nz+q_{n-1}}
.$$

We now define the notation we will use to talk about symmetry in
continued fractions.  Let $\wfor$ denote the word $a_1,a_2,\dots,a_n$.
We use $\wbac$ to denote the word $a_n,a_{n-1},\dots,a_1$, and
$-\wbac$ to denote $-a_n,-a_{n-1},\dots,-a_1$.  This notation will be
used mainly in continued fractions.  (This is the notation of
\cite{vdp}, where \propref{folding} was first stated in this form.)

We now give a proof of \propref{folding}, known as the Folding Lemma,
for completeness and in preparation for the proof in \secref{genfol} of
a generalization.

\begin{prop}
\label{folding}
$$\frac{p_n}{q_n} + \frac{(-1)^n}{xq_n^2} = [a_0,\wfor,x,-\wbac] $$
\end{prop}

\begin{pf}
First, note that $$f_n^{-1}(z) = [0,-a_n,-a_{n-1},\dots,-a_0+z].$$
In particular,
$f_n^{-1}(\infty) = [0,-\wbac]$.  However, it is easy to see that
$$f_n^{-1}(z) = \frac{-q_{n-1}z+p_{n-1}}{q_nz-p_n}.$$
Hence, $[x,-\wbac] = x-q_{n-1}/q_n$, and
$$[a_0,\wfor,x,-\wbac] = \frac{p_n(x-q_{n-1}/q_n) +
p_{n-1}}{q_n(x-q_{n-1}/q_n)+q_{n-1}}
= \frac{p_nq_nx - (p_nq_{n-1}-q_np_{n-1})}{q_n^2x}.$$
The desired result follows immediately.
\end{pf}

A continued fraction $[a_0,\wfor,x_1,-\wbac]$ is said to be folded.
We generalize this notion as follows.  A folded continued fraction has
2-fold symmetry.  We say that $[a_0,\wfor,x_1,-\wbac,x_2,\wfor]$ has
3-fold symmetry, $[a_0,\wfor,x_1,-\wbac, x_2, \wfor, x_3,-\wbac]$ has
4-fold symmetry, etc.  The Folding Lemma generalizes nicely to
$k$-fold symmetry.  We deal with the generalization in
\secref{genfol}.

For a more involved discussion of the Folding Lemma, see \cite{vdp}.
Folded continued fractions were discovered independently by
Kmo\v{s}ek and Shallit.  See
\cite{kmosek,shallit1,shallit2}.  The folding lemma first appears in
\cite[p.~209]{france}, and more explicitly in \cite[p.~332]{bfrance}.

\section{Applications of the Folding Lemma}

As an example of the sort of series with which we will deal in
\secref{appl}, we discuss an example due to Shallit (see
\cite{shallit2}).  This will serve as preparation for applying similar
techniques later.  In addition, we will need the last result of this
section.

Folded continued fractions were originally studied to explain the
simple continued fractions of certain series.  For example, the
Liouville number $\sum_{n=0}^{\infty}\frac{1}{10^{n!}}$ has an amazing
continued fraction expansion: $$\sum_{n=0}^\infty\frac{1}{10^{n!}} =
[0, 4, 1, 3, 5, 99, 1, 4, 3, 1, 4, 999999999999, 1, 3, 1, 3, 4, 1, 99,
5, 3, 1, 4,\ldots]$$

The Folding Lemma applies to the partial sums, giving a continued
fraction expansion for each of them in terms of that of the one
before.  For example, we easily get that
$$\sum_{n=0}^4\frac{1}{10^{n!}} =
[0,5,-4,-5,-100,5,4,-5,-1000000000000,5,-4,-5,100,5,4,-5].$$ This
gives the nearest integer continued fraction for the series, but it
would be more interesting to get the simple continued fraction.  (In a
simple continued fraction, all of the partial quotients are integers,
and all are positive, except perhaps the first.)

In fact, we can easily get the simple continued fraction expansion of
any number with $k$-fold symmetry, assuming that $[a_0,\wfor]$ is a
simple continued fraction.  First, note that since
$[\ldots,a,0,b,\ldots] = [\ldots,a+b,\ldots]$, it is never a problem
to have 0 occur in a continued fraction.  We can apply
$[\ldots,a,-\beta] = [\ldots,a-1,1,\beta-1]$ to get rid of negatives.
(This is the method of \cite{special}.)  For example,
$[a_0,\wfor,x_1,-\wbac] =
[a_0,a_1,\dots,a_n,x_1-1,1,a_n-1,a_{n-1},\dots,a_1].$ Similarly, in
the case of 3-fold symmetry we get
$[a_0,a_1,\dots,a_n,x_1-1,1,a_n-1,a_{n-1},\dots,a_2,a_1-1,1,x_2-1,
a_1,\dots,a_n].$ Note also that changing the sign of some $x_i$ is
easily handled, since $[\ldots,a,-x-1,1,a-1,\ldots] =
[\ldots,a-1,1,x-1,a,\ldots].$  (Note that one can remove negatives
from any continued fraction, not just one with $k$-fold symmetry.)

In this way, we get the simple continued fraction for the sum.  This
explains the amazing continued fraction noted above.  Since the
continued fraction expansion of each partial sum arises from that of
the previous partial sum by an application of the Folding Lemma, we
say that the continued fraction has iterated 2-fold symmetry.
Iterated $k$-fold symmetry is defined analogously.

Note that the Folding Lemma applies not only to $\Q$, but also to
$\Q(x)$.  If $\ell$ is an integer greater than 1, then the series
$$\sum_{n=0}^\infty\frac{1}{x^{\ell^n}}$$
converges to a formal Laurent series.  The Folding Lemma shows that it
has iterated 2-fold symmetry in its continued fraction.

In addition, the partial quotients are in $\Z[x]$.  In the terminology
of \cite{special}, the continued fraction is specializable.  We can
specialize $x$ to any positive integer greater than 1, and get a
simple continued fraction with iterated 2-fold symmetry.  In a typical
continued fraction with polynomial coefficients, the partial quotients
will have non-integral coefficients, so we will not get integral
partial quotients if we attempt to specialize.  (This point of view
was first taken in \cite{products}, although these series had been
studied earlier.)

Now consider $f(x) \in \Z[x]$ with degree greater than 1,
such that $f(x) \equiv 0 \pmod{x^2}$.  If we define the iterates
$f^0(x)=x$, and $f^n(x)=f(f^{n-1}(x))$, then
the Folding Lemma implies that
$$
\sum_{n=0}^\infty\frac{1}{f^n(x)}
$$
has a specializable continued fraction with iterated 2-fold symmetry.
In \secref{specializability}, we will determine all polynomials $f(x)
\in \Z[x]$ for which this sum is specializable.

\section{Modified Continuants}

To state the Generalized Folding Lemma, we will need to use continuant
polynomials.  Ours will sometimes have to be modified, though, to have
a sign alternation.  Fix $\varepsilon = \pm 1$.  We define modified
continuants $K'$ by
\begin{enumerate}
\item $K'(x_1) = x_1.$

\item $K'(x_1,x_2) = x_1x_2+\varepsilon.$

\item For each $k \ge 2$, $$K'(x_1,x_2,\dots,x_k) = K'(x_1,x_2,\dots,
x_{k-1})x_k + \varepsilon K'(x_1,x_2,\dots,x_{k-2}).$$
\end{enumerate}
It is sometimes convenient to define $K'$ without any variables as
$K'(\,) = 1$.
Note that when $\varepsilon = 1$, these polynomials are the usual
continuants; when $\varepsilon = -1$, they are the same except for the
sign alternation. 

Because modified continuants are the same as ordinary continuants,
except for the sign alternation when $\varepsilon = -1$, one would
expect the theory of modified continuants to be nearly the same as
that of continuants, and indeed it is.  For the theory of continuants,
see \cite{concrete}.  All of our results for modified continuants are
based directly on analogous results for continuants given in
\cite{concrete}.

One fundamental fact about modified continuants is that
$K'(x_1,x_2,\dots,x_k)$ is the sum of the terms which can be obtained
by replacing each of a collection of disjoint pairs of consecutive
variables in $x_1x_2\ldots x_k$ with $\varepsilon$.  This is easily seen
from the definition.  From this fact, it follows that the modified
continuants also satisfy the recurrence relation
\begin{equation}
\label{secondrel}
K'(x_1,x_2,\dots,x_k) = x_1 K'(x_2,x_3,\dots,x_k) +
\varepsilon K'(x_3,x_4,\dots,x_k).
\end{equation}

\begin{lem}
\label{fraclem}
We have the continued fraction expansion
$$\frac{K'(x_1,x_2,\dots,x_m)}{K'(x_2,x_3,\dots,x_m)}
= x_1+
  {\varepsilon \over\displaystyle x_2+
   {\strut \varepsilon \over\displaystyle x_3+
{\strut \varepsilon  \over\displaystyle{\cdots+
    {\strut \varepsilon  \over\displaystyle x_{m-1} +
     {\strut \varepsilon \over x_m}}}}}}
.$$
\end{lem}

\begin{pf}
The cases in which $m \le 3$ are trivial.  Now, we continue by
induction on $m$.  The continued fraction on the right is equal to
$$\frac{K'(x_1,\dots,x_{m-1}+\varepsilon
/x_m)}{K'(x_2,\dots,x_{m-1}+\varepsilon /x_m)}.$$ 
Applying the recurrence defining $K'$ twice gives
$$\frac{K'(x_1,\dots,x_{m-2})(x_{m-1}+{\varepsilon }/{x_m}) +
\varepsilon
K'(x_1,\dots,x_{m-3})}{K'(x_2,\dots,x_{m-2})(x_{m-1}+{\varepsilon
}/{x_m}) 
+ \varepsilon K'(x_2,\dots,x_{m-3})},
$$
and then
$$\frac{K'(x_1,\dots,x_{m-1}) + \varepsilon K'(x_1,\dots,x_{m-2})/x_m}
{K'(x_2,\dots,x_{m-1}) + \varepsilon K'(x_2,\dots,x_{m-2})/x_m}.$$
Multiplying the numerator and denominator by $x_m$ and applying
the recurrence again proves the lemma.
\end{pf}

\section{The Generalized Folding Lemma}
\label{genfol}

Define $S_2 = [a_0,\wfor,x_1,-\wbac]$, $S_3 =
[a_0,\wfor,x_1,-\wbac,x_2,\wfor]$, etc.  Thus, $S_k$ is the general
form of a continued fraction with $k$-fold symmetry.  Also, from now
on set $\varepsilon = (-1)^n$ (for use in modified continuants),
and for convenience set $p = p_n$ and $q = q_n$.  Our principal result
about $k$-fold symmetry is the following theorem, the Generalized Folding
Lemma: 

\begin{thm}
\label{maintheorem}
For all $k \ge 2$,
$$S_k = \frac{p}{q} + \frac{(-1)^n K'(x_2q,x_3q,\dots,x_{k-1}q)}
{qK'(x_1q,x_2q,\dots,x_{k-1}q)}.$$
\end{thm}

\begin{pf}
We prove this by induction.  It is easy to check that the case $k=2$ is
the Folding Lemma.  Now suppose that it holds for $k-1$.

In the continued fraction by which $S_k$ is defined, if one omits the
initial $a_0,\wfor$, then one is left with a continued fraction with
$(k-1)$-fold symmetry.  Call it $S'_k$.  As in the proof of the
Folding Lemma, we know that $[x_1,-\wbac] = x_1-q_{n-1}/q_n$.  It
follows by induction that $$S'_k = x_1 - \frac{q_{n-1}}{q_n} +
\frac{(-1)^n K'(x_3q_n,\dots,x_{k-1}q_n)}{q_nK'(x_2q_n,\dots,x_{k-1}q_n)}.$$
Now, since $S_k = [a_0,\wfor,S'_k]$, we see that
$$S_k = \frac{p_n\left(x_1 - \frac{q_{n-1}}{q_n} + 
\frac{(-1)^n K'(x_3q_n,\dots,x_{k-1}q_n)}
{q_nK'(x_2q_n,\dots,x_{k-1}q_n)} \right) + p_{n-1}}
{q_n\left(x_1 - \frac{q_{n-1}}{q_n} + 
\frac{(-1)^n K'(x_3q_n,\dots,x_{k-1}q_n)}
{q_nK'(x_2q_n,\dots,x_{k-1}q_n)} \right) + q_{n-1}}.$$
This is equal to (after applying $p_nq_{n-1} - q_np_{n-1} =
(-1)^{n-1}$) $$\frac{pqx_1K'(x_2q,\dots,x_{k-1}q) +
(-1)^npK'(x_3q,\dots,x_{k-1}q) + (-1)^nK'(x_2q,\dots,x_{k-1}q)}
{q^2x_1K'(x_2q,\dots,x_{k-1}q)+(-1)^nqK'(x_3q,\dots,x_{k-1}q)},$$
which simplifies to $$\frac{p}{q} +
\frac{(-1)^nK'(x_2q,\dots,x_{k-1}q)}
{q^2x_1K'(x_2q,\dots,x_{k-1}q)+(-1)^nqK'(x_3q,\dots,x_{k-1}q)}.$$
Now applying (\ref{secondrel}) to the denominator yields the final
result.
\end{pf}

If we apply \lemref{fraclem} to $S_k$, we arrive at the following
corollary to \thmref{maintheorem}:

\begin{cor}
\label{corfrac}
We have the continued fraction expansion
$$S_k
= \frac{p}{q} + \frac{1}{q}\,\,{\strut (-1)^n\over\displaystyle x_1q+
   {\strut (-1)^n\over\displaystyle x_2q+
{\strut (-1)^n \over\displaystyle{\cdots+
    {\strut (-1)^n \over\displaystyle x_{k-2}q +
     {\strut (-1)^n\over x_{k-1}q}}}}}}.
$$
\end{cor}

\section{Further Results on $k$-fold Symmetry}

A few simple algebraic manipulations give the following equivalent
formula for $S_k$:
$$
S_k = \frac{p}{q} + \frac{(-1)^{n}}{x_1q^2}\,\, {\strut 1
\over\displaystyle {1+\frac{(-1)^{n}}{x_1x_2q^2}\,\, {\strut 1
\over\displaystyle {1+\frac{(-1)^{n}}{x_2x_3q^2}\,\, {\strut 1
\over\displaystyle {\cdots+\frac{(-1)^n}{x_{k-2}x_{k-1}q^2}}}}}}}.
$$

This expansion seems messier and less natural than that given by
\corref{corfrac}, but it has one advantage.  Set
$$X_k =
\frac{(-1)^nK'(x_3q,x_4q,\dots,x_{k-1}q)}{qK'(x_2q,x_3q,\dots,x_{k-1}q)}.$$
The expansion above makes it clear that
$$S_k = \frac{p}{q} + \frac{(-1)^n}{x_1q^2}\sum_{i=0}^\infty \left(
-\frac{X_k}{x_1}\right)^i.$$ 

One might hope that the partial sums would have interesting continued
fraction expansions.  In the case of 3-fold symmetry, there is a
simple and useful description of the continued fractions for the
partial sums.  We have $$X_3 = \frac{(-1)^nK'(\,)}{qK'(x_2q)} =
\frac{(-1)^n}{x_2q^2}.$$ Note that $1/X_3$ is an integer (assuming
that $x_2$ is an integer).  The partial sums have 4-fold symmetry, as
we see from the following proposition:

\begin{prop}
\label{parsums}
For all $j \ge 1$, if $$x_3 =
-x_1\left(\frac{1-(-x_1/X_3)^j}{1-(-x_1/X_3)}\right),$$ then 
$$\frac{p}{q}+
\frac{(-1)^n}{x_1q^2}\sum_{i=0}^j \left( -\frac{X_3}{x_1}\right)^i = 
\frac{(-1)^nK'(x_2q,x_3q)}{qK'(x_1q,x_2q,x_3q)}.$$
\end{prop}

\begin{pf}
We can and shall solve the equation for the value of $x_3$ that makes
it true.  The important point is that $x_3$ turns out to be an
integer, so 4-fold symmetry occurs. 

We now solve the equation.  It is equivalent to
$$\sum_{i=0}^j \left( -\frac{X_3}{x_1}\right)^i =
{\strut 1 \over\displaystyle {
1+\frac{(-1)^{n}}{x_1x_2q^2}\,\, {\strut 1 \over\displaystyle {
1+\frac{(-1)^{n}}{x_2x_3q^2}}}}}
.$$
Equivalently,
$$
\frac{-\frac{X_3}{x_1}
-\left(-\frac{X_3}{x_1}\right)^{j+1}}{\left(-\frac{X_3}{x_1}\right)^{j+1}-1}
=  
\frac{(-1)^{n}}{x_1x_2q^2}\,\, {\strut 1 \over\displaystyle {
1+\frac{(-1)^{n}}{x_2x_3q^2}}}
.$$
This is the same as
$$
\frac{-1+\left(-\frac{X_3}{x_1}\right)^j}{\left(-\frac{X_3}{x_1}\right)^{j+1}-1}
= 
{\strut 1 \over\displaystyle {1+\frac{(-1)^{n}}{x_2x_3q^2}}}
.$$
Finally, this is equivalent to
$$
\frac{\left(-\frac{X_3}{x_1}\right)^j-1}{\left(-\frac{X_3}{x_1}\right)^{j+1}
-\left(-\frac{X_3}{x_1}\right)^j}
=
\frac{x_2x_3q^2}{(-1)^n}
.$$
The formula for $x_3$ follows immediately.
\end{pf}

This proposition may seem complicated and perhaps uninteresting.
However, in \secref{appl} we will use it to explain an interesting
continued fraction expansion. 

\section{Applications of the Generalized Folding Lemma}
\label{appl}

We now use our results to determine some explicit continued fractions.
Note that setting $k=3$ in \thmref{maintheorem}, and possibly changing
signs, shows that $$\frac{p}{q} \pm \frac{a}{abq^2\pm1}$$ has 3-fold
symmetry for all $a$ and $b$, and all combinations of signs.  Of
course, this holds true not only in the rationals, but also in fields
of rational functions, for example.  As before, we will work in
$\Q(x)$ (and also its completion $\Q((x))$, the field of formal
Laurent series).

Now suppose that $f(x)$ is polynomial over $\Z$ of degree greater than
1 such that $f(x) \equiv 1 \pmod{x^2(x-1)}$ in $\Z[x]$.  As before,
define the iterates $f^0(x)=x$, $f^n(x)=f(f^{n-1}(x))$.  Then a
trivial induction implies that
$$f^n(x) \equiv 1 \pmod{(f^0(x)f^1(x)\ldots f^{n-1}(x))^2}.$$
Using this, we see that the series
$$\sum_{n=0}^\infty\frac{1}{f^n(x)}$$ (which converges to a formal
Laurent series in $x$) has iterated 3-fold symmetry in its continued
fraction expansion, with partial quotients in $\Z[x]$.  (Similarly,
there is also 3-fold symmetry when $f(x) \equiv -1 \pmod{x^2(x+1)}$.)

In particular, if $\ell$ is a non-zero multiple of 4, then the
Chebyshev polynomial $T_\ell(x)$ (defined by $T_\ell(x) = \cos(\ell
\cos^{-1}(x))$) satisfies the congruence condition (and is of degree
greater than 1).  Since for all $a$ and $b$, $T_{ab}(x)=T_a(T_b(x)),$
it follows that
\begin{equation}
\label{chebseries}
\sum_{n=0}^\infty\frac{1}{T_{\ell^n}(x)}
\end{equation}
has iterated 3-fold symmetry.  Also, the partial quotients have
integer coefficients, so the continued fraction is specializable.
When we specialize $x$ to any positive integer greater than 1, we get
a simple continued fraction with iterated 3-fold symmetry.  For
example,
\begin{eqnarray*}
\sum_{n=0}^\infty\frac{1}{T_{4^n}(2)} &=& [0, 1, 1, 23, 1, 2, 1,
18815, 3, 1, 23, 3, 1, 23, 1, 2, 1, 106597754640383, \\ 
&& 3, 1, 23, 1, 3, 23, 1, 3, 18815, 1, 2, 1, 23, 3, 1, 23, 1, 2, 1, 
18815, 3, 1, 23, \\
\phantom{\sum_{n=0}^\infty\frac{1}{T_{4^n}(2)}} && 3, 1, 23, 1, 2, 1,
\ldots ].
\end{eqnarray*}

When $\ell$ is not divisible by 4, the series (\ref{chebseries}) displays
more complicated behavior.  Here is an apparently
typical example of the case when $\ell$ is odd:
\begin{eqnarray*}
\sum_{n=0}^\infty\frac{1}{T_{3^n}(2)} &=& 
[0, 1, 1, 5, 1, 414, 1, 2, 4, 280903, 1, 3, 3, 207, 2, 5, 1, \\
&& 22165307996832415, 6, 2, 207, 3, 4, 140451, 1, 3, 3,\\
\phantom{\sum_{n=0}^\infty\frac{1}{T_{3^n}(2)}}&& 118, 2, 2, 1, 1, 7,\ldots]
\end{eqnarray*}
There is partial symmetry, but it breaks down.  It would be
interesting to have an explanation of this behavior, and of that for
$\ell \equiv 2 \pmod{4}.$ (By \thmref{done}, these series never have
specializable continued fractions.)

We can also use \propref{parsums} to prove that certain series have
4-fold symmetry in their continued fractions.  For example,
\propref{parsums} implies that for each $\ell$,
$$\sum_{n=0}^{3\ell+2}\frac{1}{x^{2^n}}+\sum_{n=0}^\ell\frac{1}{x^{3
\cdot 8^n}}$$ 
has iterated 4-fold symmetry.  To see this, write the sum as
$$
\sum_{n=0}^{\ell}\left(\frac{1}{x^{8^n}}+\frac{1}{x^{2\cdot 8^n}}+
\frac{1}{x^{3\cdot 8^n}}+\frac{1}{x^{4 \cdot 8^n}}\right)
$$
and apply \propref{parsums} with $q= x^{4\cdot 8^{\ell-1}}$, $j=3$,
$x_1=(-1)^n$, and $x_2 = -1$.
This explains the symmetry observed in
\begin{eqnarray*}
\label{twosum}
\sum_{n=0}^\infty\frac{1}{2^{2^n}}+\sum_{n=0}^\infty\frac{1}{8^{8^n}} &=&
[0, 1, 16, 14, 16, 1, 65792, 15, 17, 65792, 1, 16, 14, 16, 18, 14, 16 \\ 
&& 1, 65792, 16, 340282366920938463481821351505477763072, \\
\phantom{\sum_{n=0}^\infty\frac{1}{2^{2^n}}+\sum_{n=0}^\infty\frac{1}{8^{8^n}}}
&& 1, 15, 65792, 1, 16, 14, 16, \ldots].
\end{eqnarray*}
Also, one checks easily that the extremely large partial quotients are
exactly the numbers $2^{2^{3n+1}}+2^{2^{3n}}.$ 

One can generate other examples of this phenomenon.  For example,
\propref{parsums} implies that 
$$\sum_{n=0}^\infty\frac{1}{x^{6^n}}
+\sum_{n=0}^\infty\frac{1}{x^{2\cdot 6^n}} +
\sum_{n=0}^\infty\frac{1}{x^{3 \cdot 6^n}}$$ 
has iterated 4-fold symmetry.  It is unclear whether it is essentially
a coincidence that these sums have iterated 4-fold symmetry, or
whether these is some more general result along these lines.

\section{Specializability}
\label{specializability}

These results bring up some interesting related questions.  Suppose
that $f(x)$ is a polynomial over $\Z$, of degree greater than 1.
Under what conditions does
\begin{equation}
\label{itersums}
\sum_{n=0}^\infty\frac{1}{f^n(x)}
\end{equation}
have a specializable continued fraction?  Several of the sums dealt
with earlier in this paper were of this form, with the allowable
functions characterized by congruence conditions.  In this section, we
will show that in general, such a series is specializable if and only
if $f(x)$ satisfies one of fourteen congruence conditions.

All of our examples so far have had symmetry.  However, this is not
the case in general.  We will prove shortly that if $f(x) \equiv -x
\pmod{x^2}$, then the sum has a specializable continued fraction,
although $k$-fold symmetry does not occur.  The proof is fairly
typical of our methods in this section.  We first guess what the
continued fraction expansion is, and then prove it by induction.
Surprisingly, this works in every case in which the continued fraction
expansion is specializable, although several cases are tricky.  More
surprisingly, we can actually rule out every other case, and thereby
arrive at a complete classification of the polynomials for which the
series (\ref{itersums}) has a specializable continued fraction.

Before we deal with the case of $f(x) \equiv -x \pmod{x^2}$, we need
to set up some notation.  Let $f(x) = g(x)x^2-x$, and then set $A_1(x)
= -g(x)$, $A_2(x) = g(f(x))x^2$, $A_3(x) = -g(f^2(x))(f(x)/x)^2$,
$A_4(x) = g(f^3(x))x^2(f^2(x)/f(x))^2$, etc.  In general,
$$A_{2\ell}(x) = g(f^{2\ell-1}(x))x^2(f^2(x)/f(x))^2(f^4(x)/f^3(x))^2
\ldots(f^{2\ell-2}(x)/f^{2\ell-3}(x))^2,$$
and $$A_{2\ell+1}(x) = -g(f^{2\ell}(x))(f(x)/x)^2(f^3(x)/f^2(x))^2
\ldots(f^{2\ell-1}(x)/f^{2\ell-2}(x))^2.$$
(The condition on $f(x)$ implies that these are polynomials.)
Then we have the following result:

\begin{prop}
Let $f(x)=g(x)x^2-x$, and $A_i(x)$ be as defined above.
Then for each $\ell$,
$$\sum_{n=0}^\ell\frac{1}{f^n(x)} = [0,x,A_1(x),A_2(x),\dots,A_\ell(x)].$$
\end{prop}

\begin{pf}
We prove this by induction.  The base case is trivial.  Now note that
$$-x^2[A_1(f(x)),A_2(f(x)),\dots,A_\ell(f(x))] =
[A_2(x),A_3(x),\dots,A_{\ell+1}(x)].$$
This, combined with the identity
$$[0,x,-g(x),-x^2X] = 1/x + [0,f(x),X]
,$$
proves the result.
\end{pf}

We get a similar result when $f(x) \equiv x^2-x+1 \pmod{x^2(x-1)^2}$.
Let $f(x) = x^2-x+1+x^2(x-1)^2g(x)$, and set $A_1(x) = -g(x)(x-1)^2$,
$B_1(x) = -x^2$, and in general
$$A_\ell(x) =
-g(f^{\ell-1}(x))\left(\frac{f^{\ell-1}(x)-1}{xf(x)\ldots
f^{\ell-2}(x)}\right)^2,$$ 
and $$B_\ell(x) = -x^2f(x)^2f^2(x)^2\ldots f^{\ell-1}(x)^2
.$$
(As in the previous example, the condition on $f(x)$ implies that
these are polynomials.)
Then we have the following result:

\begin{prop}
Let $f(x)=x^2-x+1+x^2(x-1)^2g(x)$, and $A_i(x)$ and $B_i(x)$ be as
defined above. Then for each $\ell$,
$$\sum_{n=0}^\ell\frac{1}{f^n(x)} =
[0,x,A_1(x)-1,B_1(x),A_2(x),B_2(x),\dots,A_\ell(x),B_\ell(x)].$$ 
\end{prop}

\begin{pf}
We prove this by induction.  The base case is trivial.  Now note that
$$[A_1(f(x)),B_1(f(x)),\dots,A_\ell(f(x)),B_\ell(f(x))] =
x^2[A_2(x),B_2(x),\dots,A_{\ell+1}(x),B_{\ell+1}(x)].$$ 
This, combined with the identity
$$[0,x,A_1(x)-1,B_1(x),X] = 1/x+[0,f(x),-1+x^2X],$$
proves the result.
\end{pf}

Note that in the special case where $f(x) = x^2-x+1$, we have
$$\sum_{n=0}^\infty\frac{1}{f^n(x)} = \frac{1}{x-1}.$$
Of course, this example is rather special.

We now determine exactly which polynomials $f(x)$ are such that all of
the partial sums of $\sum_{n=0}^\infty\frac{1}{f^n(x)}$ have
specializable continued fractions.  Before doing this, however, we
need to discuss a few points about proving non-specializability.

Suppose that one has a continued fraction expansion for a formal
Laurent series.  If one of the partial quotients does not have
integral coefficients, it is not necessarily the case that the Laurent
series does not have a specializable continued fraction.  To see this,
note that the identity $[a+1/b,c]=[a,b,-(c+b)/b^2]$ shows that when
the constant term of a partial quotient is not an integer, one can on
occasion adjust the continued fraction to make it specializable.  For
example, $[x-1/3,9x^2+3] = [x,-3,-x^2]$.  However, it is not hard to
see that if the first partial quotient without integer coefficients
has a non-integral coefficient other than the constant term, then the
Laurent series has no specializable continued fraction.  This
observation will suffice for all of the examples we will consider.

Also, we will look at continued fractions involving several variables.
When one wants to prove that such a continued fraction is not
specializable, one must be careful about 0 as a partial quotient.  If
one of the partial quotients becomes 0 for certain values of some of
the variables, then its two neighboring partial quotients add when the
variables assume those values.  This could make the continued fraction
specializable.  In general, it will be clear that this doesn't happen,
however.

We begin with the following easy lemma.

\begin{lem}
Let $f(x)$ be a (non-zero) polynomial over $\Z$. Then $1/x + 1/f(x)$
has a specializable continued fraction iff $f(x)$ is congruent modulo
$x^2$ to one of $$0,-1,1,-x,-x-1,-x+1,-2x,-2x-1,-2x+1.$$
\end{lem}

\begin{pf}
Suppose that $f(x) = g(x)x^2 + bx + a$.  Then
$$ \frac{1}{x}+\frac{1}{f(x)} = [0,x,-g(x),
-\frac{x}{b+1}+\frac{a}{(b+1)^2},-\frac{(b+1)^3x}{a^2}-\frac{(b+1)^2}{a}]
.$$ 
From this, we see that unless $b=0$, $b=-1$, or $b=-2$, there is
no specializable continued fraction for $1/x+1/f(x)$.  Now, we just
check each case.

When $b=0,$ $1/x+1/f(x) = [0,x,-g(x),-x+a, -x/a^2-1/a].$ When $b=-1$,
$1/x+1/f(x) = [0,x, -g(x),-x^2/a].$ When $b=-2$, $1/x+1/f(x) =
[0,x,-g(x),x+a,x/a^2-1/a].$ From these, we see that in each case, we
must have $a=-1$, $a=0$, or $a=1$ to have specializability, and that
in each of those cases, $1/x+1/f(x)$ does have a specializable
continued fraction.
\end{pf}

\begin{lem}
\label{oneproof}
Let $f(x)$ be a (non-zero) polynomial over $\Z$, and suppose that
$f(x) = 1+kx^2+x^2(x-1)g(x)$ with $g(x) \in \Z[x].$  Then
$1/x+1/f(x)+1/f^2(x)$ has a specializable continued fraction iff $k=0.$
\end{lem}

\begin{pf}
We look at the more general sum
$$\frac{1}{x} + \frac{1}{f(x)} + \frac{1}{1+kf(x)^2+f(x)^2(f(x)-1)G(x)}
.$$ (Note that we do not require that $G(x) = g(f(x)).$)

If $k=0$, this has the continued fraction expansion
$[0, x, -g(x)(x-1), -x+1, -x-1, G(x)g(x)(x-1), x+1, x-1, g(x)(x-1),
-x+1, -x-1, g(x)(x-1), x-1, x+1].$  If $k=1$, then it has the
expansion $[0, x, -g(x)(x-1)-1, -x+1, -x-1, G(x)g(x)(x-1)+G(x), x,
-g(x)x^3+(g(x)-1)x^2, -x, g(x)(x-1)+1, x-2, x/4+1/2],$ which cannot be made
specializable.  Finally, for all other values of $k$, it has the
expansion
$[0,x,-g(x)(x-1)-k,-x+1,-x-1,G(x)g(x)(x-1)+G(x)k,x+1-k,{\frac
{x}{k^{2}-2k+1}}+\frac{1}{k-1},g(x)\left (k^{2}-2k+1\right )^{2}x+\left 
(k^{2}-2k+1\right)
(k-g(x)+k^{3}+2\,kg(x)-2\,k^{2}-k^{2}g(x)), -{\frac
{x}{k^{2}-2k+1}}-\frac{1}{k-1},-x-1+k,g(x)(x-1)+k,x-1-k,{\frac
{x}{1+2\,k+k^{2}}}+\frac{1}{k+1}],$ which also cannot be made
specializable.
\end{pf}

Now note that changing $f(x)$ to $-f(-x)$ has no effect on the
specializability of the continued fractions for the partial sums we
are studying.  Thus, we immediately get the following lemma:

\begin{lem}
Let $f(x)$ be a (non-zero) polynomial over $\Z$, and suppose that
$f(x) = -1 + kx^2 +x^2(x+1)g(x)$ with $g(x) \in \Z[x].$  Then
$1/x+1/f(x)+1/f^2(x)$ has a specializable continued fraction iff
$k=0.$
\end{lem}

\begin{lem}
Let $f(x)$ be a (non-zero) polynomial over $\Z$, and suppose that
$f(x) = -2x + x^2g(x)$ with $g(x) \in \Z[x].$  Then $1/x + 1/f(x)
+1/f^2(x)$ never has a specializable continued fraction.
\end{lem}

\begin{pf}
As in the proof of \lemref{oneproof}, we generalize slightly,
and look at
$$\frac{1}{x}+\frac{1}{f(x)}+\frac{1}{-2f(x)+f(x)^2G(x)},$$
where $G(x)$ need not be equal to $g(f(x)).$  This has the continued
fraction expansion $[0, x, -g(x), x, -G(x), -x, \break -g(x), x/3]$, which
cannot be made specializable.
\end{pf}

Similar arguments, best carried out using computer algebra software,
prove the following two lemmas:

\begin{lem}
Let $f(x)$ be a (non-zero) polynomial over $\Z$, and suppose that
$f(x) \equiv -2x+1 \pmod{x^2}.$  Then $1/x+1/f(x)+1/f^2(x)$ has a
specializable continued fraction iff $f(x) \equiv -x^3+3x^2-2x+1
\pmod{x^2(x-1)^2}$ or $f(x) \equiv x^2-2x+1 \pmod{x^2(x-1)^2}.$
\end{lem}

\begin{lem}
Let $f(x)$ be a (non-zero) polynomial over $\Z$, and suppose that
$f(x) \equiv -x+1 \pmod{x^2}.$  Then $1/x+1/f(x)+1/f^2(x)$ has a
specializable continued fraction iff $f(x) \equiv x^2-x+1, x^3-x+1,
x^3-x^2-x+1, {\text{or}}\,\, {-x^3+2x^2-x+1} \pmod{x^2(x-1)^2}.$
\end{lem}

(Of course, changing $f(x)$ to $-f(-x)$ gives analogous results for
the cases in which $f(x) \equiv -2x-1 \pmod{x^2}$ or $f(x) \equiv -x-1
\pmod{x^2}.$)

In fact, when $f(x) \equiv x^3-x+1 \pmod{x^2(x-1)^2},$ one can check that
$1/x + 1/f(x) +1/f^2(x) + 1/f^3(x)$ does not have a specializable continued
fraction.  This, combined with the preceding lemmas, proves the necessity of
the conditions of the following theorem:

\begin{thm}
\label{big}
Let $f(x)$ be a polynomial over $\Z$, of degree greater than 1.  All
of the partial sums of $\sum_{n=0}^\infty\frac{1}{f^n(x)}$ have
specializable continued fractions if and only if $f(x)$ satisfies
one of the following congruences:
\begin{enumerate}
\item $f(x) \equiv 0 \pmod{x^2}$.
\item $f(x) \equiv -x \pmod{x^2}$.
\item $f(x) \equiv 1 \pmod{x^2(x-1)}$.
\item $f(x) \equiv -1 \pmod{x^2(x+1)}$.
\item $f(x) \equiv x^3-x^2-x+1 \pmod{x^2(x-1)^2}$.
\item $f(x) \equiv -x^3+2x^2-x+1 \pmod{x^2(x-1)^2}$.
\item $f(x) \equiv -x^3+3x^2-2x+1 \pmod{x^2(x-1)^2}$.
\item $f(x) \equiv x^3+x^2-x-1 \pmod{x^2(x+1)^2}$.
\item $f(x) \equiv -x^3-2x^2-x-1 \pmod{x^2(x+1)^2}$.
\item $f(x) \equiv -x^3-3x^2-2x-1 \pmod{x^2(x+1)^2}$.
\item $f(x) \equiv x^2-x+1 \pmod{x^2(x-1)^2}$.
\item $f(x) \equiv x^2-2x+1 \pmod{x^2(x-1)^2}$.
\item $f(x) \equiv -x^2-x-1 \pmod{x^2(x+1)^2}$.
\item $f(x) \equiv -x^2-2x-1 \pmod{x^2(x+1)^2}$.
\end{enumerate}
\end{thm}

\begin{pf}
We have already proved the necessity that one of these congruences
hold.  We now prove its sufficiency.

Because of our results so far, and the symmetry between $f(x)$ and
$-f(-x)$, we see that we need only deal with (5), (6), (7), and (12).
Set $S_n(x) = 1/x+1/f(x)+\cdots+1/f^n(x),$ and $Q_n(x) = xf(x)\cdots f^n(x).$

We begin with (6).  Suppose that $f(x)
= -x^3+2x^2-x+1 +x^2(x-1)^2g(x).$  First, define
the polynomials $A_0(x) = 1$, $A_1(x) = (1-f(x))/x$, and for $n \ge 2$
$$ A_n(x) = \frac{1-f^n(x)}{xf(x)\ldots f^{n-1}(x)} +
f^n(x)f^{n-1}(x)A_{n-2}(x).$$ 
(Because $f(x)-1 = x(x-1)^2(xg(x)-1),$
these are indeed polynomials.  In fact, we see from this that for $n >
1,$ $f^n(x)-1$ is divisible by $Q_{n-1}(x)Q_{n-2}(x)$.)

In fact, $A_n(x)$ is the denominator of the penultimate convergent of
$S_n(x)$.  We will see later that the continued fraction expansion of
$S_n(x)$ has even length.
If $S_n(x) = p/q$ with
$p$ and $q$ coprime, then $p'/q'$ is the penultimate convergent iff
$pq'-qp' = 1,$ and also $\deg q' < \deg q$ and $\deg p' < \deg
p$. In this case, it is easiest to check simply that
$q'(p/q)-1/q$ is a polynomial of degree less than $\deg p =
\deg q -1.$  Of course, the denominator of $S_n(x)$ is easily
seen to be $Q_n(x)$.  (This is
the method used in \cite{special}.)

Our proof that $A_n(x)$ is the denominator of the penultimate
convergent will be by induction.  The first two cases are easy to
check.  Also, note that if $S_n(x)A_n(x)-1/Q_n(x)$ is a polynomial,
then it automatically has the right degree.  Now, we express it as
$$
\left(S_{n-2}(x)+\frac{1}{f^{n-1}(x)}+
\frac{1}{f^n(x)}\right)\left(\frac{1-f^n(x)}{Q_{n-1}(x)}
+f^n(x)f^{n-1}(x)A_{n-2}(x)\right) -\frac{1}{Q_n(x)}.$$ 
The product $(1/f^{n-1}(x)+1/f^n(x))f^n(x)f^{n-1}(x)A_{n-2}(x)$ is a
polynomial, as is $S_{n-2}(x)(1-f^n(x))/Q_{n-1}(x)$.  (The latter one
is a polynomial since $f^n(x)-1$ is divisible by
$Q_{n-1}(x)Q_{n-2}(x)$.)  Also, by induction,
$S_{n-2}(x)A_{n-2}(x)f^n(x)f^{n-1}(x)$ is a polynomial plus
$f^n(x)f^{n-1}(x)/Q_{n-2}(x)$.

Thus, we need only show that
\begin{equation}
\label{messyind}
\frac{f^n(x)f^{n-1}(x)}{Q_{n-2}(x)} +
\left(\frac{1}{f^{n-1}(x)}+
\frac{1}{f^n(x)}\right)\left(\frac{1-f^n(x)}{Q_{n-1}(x)}\right)
- \frac{1}{Q_n(x)}
\end{equation}
is a polynomial, or equivalently that
$$ f^n(x)f^{n-1}(x) + \left(\frac{1}{f^{n-1}(x)}+\frac{1}{f^n(x)}\right)
\left(\frac{1-f^n(x)}{f^{n-1}(x)}\right)
-\frac{1}{f^n(x)f^{n-1}(x)}
$$
is a polynomial which is divisible by $Q_{n-2}(x)$.  This follows
immediately from that fact that when $n=2$, it is a polynomial which
is divisible by $f(x)-1$, because one sees then that it is always a
polynomial divisible by $f^{n-1}(x)-1$, which is divisible by
$Q_{n-2}(x)Q_{n-3}(x)$.

Now let
$X_1(x) = -g(x)(x-1)^2+x-2$, and for $n \ge 2$,
$$X_n(x) =
\frac{A_{n-2}(x)}{Q_{n-2}(x)}-\frac{A_{n-1}(x)}{Q_{n-1}(x)}+
\frac{1-f^n(x)}{Q_{n-1}(x)^2}.$$
(One can prove by induction that this is a polynomial, by looking
at the differences $X_n(x)-X_{n-2}(x)$ and showing that they
are polynomials in a way like that used to prove that (\ref{messyind})
is a polynomial.)

Note that $S_1(x) = [0,x,X_1(x),-x^2].$ Now suppose that $\sfor$ is
such that for each $n$, $S_n(x) = [0,\sfor].$ Then for $n \ge 2$ we
have $S_n(x) = [0,\safor,X_n(x),-\sbbac,f^{n-1}(x)^2,\sbfor].$
(Strictly speaking, one should probably define $\sfor$ inductively as
$\safor,X_n(x),-\sbbac,f^{n-1}(x)^2$, $\sbfor$, and then note that
$S_n(x) = [0,\sfor].$)

This is not hard to prove.  Of course, 
$$[X_n(x),-\sbbac] = X_n(x)-A_{n-2}(x)/Q_{n-2}(x).$$
It follows from the Folding Lemma that
$$[X_n(x),-\sbbac,f^{n-1}(x)^2,\sbfor] = X_n(x) -
\frac{A_{n-2}(x)}{Q_{n-2}(x)}-\frac{1}{f^{n-1}(x)^2Q_{n-2}(x)^2}
.$$
From this, one sees that if one subtracts $S_{n-1}(x)$ from the whole
continued fraction, one gets
$$
\frac{1}{Q_{n-1}(x)\left(-Q_{n-1}(x)\left(X_n(x)-{\displaystyle
A_{n-2}(x) \over \displaystyle Q_{n-2}(x)}-\frac{\displaystyle
1}{\displaystyle f^{n-1}(x)^2Q_{n-2}(x)^2}\right)-A_{n-1}(x)\right)}.
$$
If one applies the definitions of $X_n(x)$ and of $Q_n(x)$, this
simplifies to $1/f^n(x)$.  (There is no need to use the recurrence
defining $A_\ell(x)$ in these manipulations, because all of the
occurrences of $A_\ell(x)$ (for various $\ell$) cancel.)  We omit the
details.  This proves the sufficiency of case (6).

We now deal with case (12).  Suppose that $f(x) = (x-1)^2(1+g(x)x^2).$
First, define the polynomials (for $n \ge 1$)
$$
A_n(x)  = \frac{f^n(x)-1}{f^{n-1}(x)} + f^n(x)f^{n-1}(x)S_{n-1}(x)
.$$
(One can check easily that these are polynomials.)

As before, $A_n(x)$ is the denominator of the penultimate convergent
of $S_n(x)$.  We will prove that the continued fraction expansion of
$S_n(x)$ has odd length.  Assuming this, it is easy to prove that
$A_n(x)$ is indeed this denominator.
Note that the denominator of $S_n(x)$ is $f^n(x)f^{n-1}(x)$.
Let $B_n(x) = A_n(x)S_n(x)+1/(f^n(x)f^{n-1}(x)).$  We need to show that
$B_n(x)$ is a polynomial.  To do this by induction, we look at the difference
$$B_n(x)-B_{n-1}(f(x)) = \left(A_{n-1}(f(x))+\frac{f^n(x)f^{n-1}(x)}{x}\right)
\left(S_{n-1}(f(x))+\frac{1}{x}\right) - A_{n-1}(f(x))S_{n-1}(f(x)),
$$
and show that it is a polynomial.  It is easy to check that
$f^n(x)f^{n-1}(x)/x^2$ is a polynomial.  Thus, this reduces to
checking that $A_{n-1}(f(x))+S_{n-1}(f(x))f^n(x)f^{n-1}(x)$ is a
polynomial divisible by $x$.  Call this $C_n(x)$.  Then $C_{n+1}(x) =
C_{n}(f(x))+f^n(x)f^{n-1}(x)/x+f^{n+1}(x)f^n(x)/f(x),$ and one sees by
induction that these are polynomials divisible by $x$.

Now, for $n \ge 2$, define
$$
X_n(x) = g(f^{n-1}(x))\left(\frac{f^{n-1}(x)-1}{f^{n-2}(x)}\right)^2
$$
and for $n \ge 3$,
$$
Y_n(x) = -\frac{f^{n-1}(x)}{f^{n-3}(x)^2}
.$$
Also define $Y_2(x) = f(x)$.

Now note that $S_1(x) = [0,x,-g(x)(x-1)^2-1,x+1,x-1]$ and $S_2(x) = 
[0, x, -g(x)(x-1)^2-1, x+1, x-1, X_2(x), -x, f(x)+1, x, -g(x)(x-1)^2-1, x]
$.  Suppose that for each $n$, $S_n(x) = [0,\sfor]$.  Then for $n \ge
3$,  we have
\begin{equation}
\label{cfracexp}
S_n(x) =
[0,\safor,X_n(x),-\sbfor,Y_n(x),\scfor,Y_{n-1}(x),-\scbac,X_{n-1}(x),\sbbac].
\end{equation}

This can be proved completely straightforwardly.  The case $n=3$ seems
to need to be treated differently, since it involves $Y_2(x)$, which
does not quite fit the pattern defining $Y_n(x)$ in general, but the
same methods suffice to prove it.

In the general case, the techniques we have been using show that the
continued fraction in (\ref{cfracexp}) is equal to a certain rational
function of iterates of $f(x)$, $g(x)$ applied to iterates of
$f(x)$,  and partial sums of the series.  (Note that the numerators
and denominators of the last two convergents of $S_n(x)$ can be
expressed in terms of these.)

When one expresses each of the partial sums in terms of $S_n(x)$ (and
the iterates of $f(x)$) and subtracts off $S_n(x)$, one is left with a
rational function of the iterates of $f(x)$, and $g(x)$ applied to
them.  (The partial sums all disappear by cancellation.)  One can than
express each of the iterates of $f(x)$ in terms of the first one to
appear.  This gives a rational function of one iterate of $f(x)$, as
well as $g(x)$ applied to several iterates.  This rational function
vanishes identically.  (This would be tedious to check by hand, but is not
very difficult to check using computer algebra software.)

Now, we deal with case (5).  Suppose that $f(x) =
(x-1)^2(1+x+x^2g(x)).$  Now define $A_1(x) = (f(x)-1)/x$, and
for $n \ge 2,$
$$
A_n(x) = \frac{f^n(x)-1}{f^{n-1}(x)}+f^n(x)f^{n-1}(x)S_{n-2}(x)
.$$
One can check, as before, that this is the denominator of the
penultimate convergent of $S_n(x)$.  We will prove shortly that the
continued fraction has odd length for $n \ge 2$.  Note that the
denominator of $S_n(x)$ is $f^n(x)f^{n-1}(x).$

Now define $X_n(x)$ for $n \ge 3$ to be
$$
g(f^{n-1}(x))\left( \frac{f^{n-1}(x)-1}{f^{n-2}(x)}\right)^2 +
g(f^{n-2}(x))(f^{n-2}(x)-1)^2 + f^{n-2}(x)-1
,$$
$Y_n(x)$ for $n \ge 4$ to be
$$
\left( \frac{f^{n-2}(x)-1}{f^{n-3}(x)}\right)^2\!
(g(f^{n-2}(x))(f^{n-1}(x)+(f^{n-2}(x)+1)(f^{n-2}(x)-1)^2) + f^{n-2}(x)^2-2)
,$$
and $Z_n(x)$ for $n \ge 6$ to be
$$
-Y_{n-2}(x) - \frac{f^{n-2}(x)f^{n-3}(x)^2}{f^{n-4}(x)^2f^{n-5}(x)^2}
.$$
(These are easily seen to be polynomials.)

Also define $X_1(x) = -g(x)(x-1)^2-x+1$, and $X_2(x)$ following the definition
above, except multiplied by $-1$.  Define $Y_1(x) = -x^2$, $Y_2(x) = 
-g(x)(x-1)^4(g(x)x^2+2x+2)-x^4+2x^3+x^2-4x+1$, and $Y_3(x)$ following
the definition above, except multiplied by $-1$.  Finally, define
$Z_3(x) = x^2(f(x)+1)$, $Z_4(x) = (f(x)^2(f(f(x))+1)+2x-1)/x^2$, and
$Z_5(x)$ following the
definition above, except with the second term multiplied by $-1$.

Now note that $S_1(x) = [0,x,X_1(x),Y_1(x)]$ and $S_2(x) = [0,x,X_1(x),Y_1(x),
X_2(x), x-1,x+1,Y_2(x),x]$.
Suppose that for each $n$, $S_n(x) = [0,\sfor].$  Then for $n \ge
3$ we have
$$S_n(x) = [0,\safor,X_n(x),-\sbfor,0,\sdfor,Y_{n-2}(x),0,Z_n(x),
-\sdbac,Y_n(x),\sbbac].$$
(When $n=3$, $\smfor$ and $-\smbac$ appear here.  They should be
interpreted to be empty words. Also, note that because $0$ appears in
this expansion, the continued fraction collapses somewhat.
However, it remains specializable.)

As before, this is straightforward to prove.  The methods used for the
proof of the previous cases also work in this case.  (It is best to
use computer algebra software for the calculations.)

Finally, we deal with case (7).
First note that given any continued fraction $[a_0,a_1,\dots,a_n]$
with convergents $p_i/q_i$, and given any $X$, we have
\begin{equation}
\label{help}
[a_0,\wfor,X,\wbac,1,-\wfor] = \frac{p_n}{q_n} +
\frac{1}{q_n(Xq_n+2q_{n-1})(-1)^n+1}. 
\end{equation}
(This is easy to prove using the techniques used to prove the
Generalized Folding Lemma.) 

Now, suppose that $f(x) = -x^3+3x^2-2x+1+g(x)x^2(x-1)^2$.  Note that
$f(x)-1 = x(x-1)(g(x)(x^2-x)-x+2),$ and that the denominator of
$S_n(x)$ is $Q_n(x)$.  We will show that that continued fraction of
$S_{n+1}(x)$ arises from that of $S_n(x)$ through application of
(\ref{help}), for a suitable choice of $X$.

To prove this, we use the fact that for any continued fraction,
$p_kq_{k-1} \equiv (-1)^{k-1} \pmod{q_k}.$ (This follows from the
equation $p_kq_{k-1}-q_kp_{k-1}=(-1)^{k-1}$.)  From this, one sees
immediately that (\ref{help}) is applicable (with $X$ a polynomial) if
and only if
$$
S_n(x)(f^{n+1}(x)-1) \equiv -2 \pmod{Q_n(x)}
.$$
(It is easy to check that $Q_n(x)$ divides $f^{n+1}(x)-1$, so
$S_n(x)(f^{n+1}(x)-1)$ is a polynomial.) 

However, this is easy to prove by induction.  The base case is
trivial.  Now, assume that it holds for $n-1$.  Using induction and
the fact that $f^n(x)-1$ divides $f^{n+1}(x)-1$, we have
$$S_{n-1}(x)(f^{n+1}(x)-1) \equiv
-2\left(\frac{f^{n+1}(x)-1}{f^n(x)-1}\right) \pmod{Q_n(x)}.
$$ 
Hence,
$$
S_n(x)(f^{n+1}(x)-1) \equiv -2\left(\frac{f^{n+1}(x)-1}{f^n(x)-1}\right) +
\frac{f^{n+1}(x)-1}{f^n(x)}
\pmod{Q_n(x)}.
$$
This simplifies to $-f^n(x)(f^n(x) - 1)((f^n(x)+1)g(f^n(x)) - 1) - 2$.
Because $Q_{n-1}(x)$ divides $f^n(x)-1$, this is $-2$ modulo $Q_n(x)$,
as desired.  Hence, the continued fraction of $S_{n+1}(x)$ is
determined by that of $S_n(x)$ by (\ref{help}).  Since the continued
fraction of $S_1(x)$ is $[0,x,-g(x)x^2+(2g(x)+1)x-g(x)-3,x+1,x-1]$,
all of these continued fractions are specializable.

This completes the proof.
\end{pf}

\begin{lem}
Let $f(x)$ be a quadratic polynomial over $\Z$.  Then
$\sum_{n=0}^\infty\frac{1}{f^n(x)}$ has a specializable continued
fraction iff $f(x) = kx^2$ for some $k$, $f(x) = kx^2-x$ for some $k$,
or $f(x)$ is one of $x^2-x+1$, $x^2-2x+1$, $-x^2-x-1$, or $-x^2-2x-1$.
\end{lem}

(This lemma is proved the same way as the necessity of the conditions
in \thmref{big}.  We omit the details.  Essentially, one looks at
partial sums and shows that they do not have specializable continued
fractions, unless $f(x)$ belongs to one of the six cases of the lemma.
Arguments like those used to prove the next lemma show that the
beginnings of the continued fractions for these partial sums coincide
with the beginning of the continued fraction for the entire series.
In each case, the non-specializability of the continued fraction for
the partial sum occurs near enough to its beginning to imply that the
entire series cannot have a specializable continued fraction
expansion.)

\begin{lem}
Let $f(x)$ be a polynomial over $\Z$, of degree greater than 2.  Then
the partial sums of $\sum_{n=0}^\infty\frac{1}{f^n(x)}$ are
convergents to its continued fraction.
\end{lem}

\begin{pf}
Suppose that $\deg f(x) = \ell$.  Since denominator of the partial
sum $\sum_{n=0}^k{1/f^n(x)}$ clearly divides $f^0(x)f^1(x)\ldots f^k(x)$, it
has degree at most $(\ell^{k+1}-1)/(\ell-1)$.  Since $\ell \ge
3$, this is strictly less than $\frac{1}{2}\deg f^{k+1}(x)$.

Now note that if we define a valuation on $\Q(x)$
by $|a(x)/b(x)| = 2^{\deg a(x) - \deg b(x)},$ 
then it is easy to see that for $a,\varepsilon \in \Q(x)$,
$|a\varepsilon| < 1$ implies
$\varepsilon + 1/a$ = $1/(a+\varepsilon')$ with $|\varepsilon'| =
|\varepsilon a^2|$.  From this, it follows immediately that
$\sum_{n=0}^k{1/f^n(x)}$ is a convergent of
$\sum_{n=0}^{k+1}{1/f^n(x)}$.  This proves our lemma.
\end{pf}

Combining the last three results yields the following theorem:

\begin{thm}
\label{done}
Let $f(x)$ be a polynomial over $\Z$, of degree greater than 1.  Then
$\sum_{n=0}^\infty\frac{1}{f^n(x)}$ has a specializable continued
fraction if and only if $f(x)$ satisfies one of the fourteen congruences
listed in the statement of \thmref{big}.
\end{thm}

Another interesting example is $f(x) = x^2+x-1$.  The continued
fraction in this case is not specializable, but no denominators
greater than 2 occur.  Doubling the continued fraction eliminates all
but one of them.  If a continued fraction is not specializable, but
the coefficients of its partial quotients have denominators at most 2,
then we call it semi-specializable.  Semi-specializable continued
fractions are almost as interesting as specializable ones, because it
is fairly straightforward to transform a continued fraction whose
partial quotients are halves of integers to a simple continued
fraction.  It might be interesting to find an analogue of
\thmref{done} for semi-specializability, because it is conceivable
that there would be a simpler or more interesting characterization.

\section{Other Forms of Folding Symmetry}
\label{other}

In addition to folding symmetry, at least one other similar form of
symmetry occurs in some continued fractions, this time in products
rather than sums.  Define $\afor$ to be the word $a_0,a_1,\dots,a_n$,
and define $\abac$ to be $a_n,a_{n-1},\dots,a_0$.  Thus, $\afor =
a_0,\wfor$.  This notation will be more convenient than the earlier
notation.

Note that since $[\,\abac] = p_n/p_{n-1}$, $\afor = \abac$ if and only
if $p_{n-1} = q_n$.  (To compute $[\,\abac]$, note that $[\,\abac] =
-1/f_n^{-1}(0) = p_n/p_{n-1}.$) We now use this to prove the following
proposition:

\begin{prop}
\label{strange}
If $\afor = \abac$ then
$$[\,\afor,x,\abac] = \frac{p_n}{q_n}\left(
1+\frac{(-1)^n}{q_n(xp_n+2q_n)+(-1)^{n-1}}\right).
$$
\end{prop}

\begin{pf}
We have
$$[\,\afor,x,\abac] =
\frac{p_n(x+p_{n-1}/p_n)+p_{n-1}}{q_n(x+p_{n-1}/p_n) + q_{n-1}}.$$
Now if we apply the relation $p_nq_{n-1} = (-1)^{n-1} + q_np_{n-1}$,
we get 
$$\frac{p_n^2x + 2p_np_{n-1}}{p_nq_nx + 2q_np_{n-1} + (-1)^{n-1}},$$
and substituting $p_{n-1} = q_n$ gives
$$\frac{p_n^2x + 2p_nq_n}{p_nq_nx + 2q_n^2 + (-1)^{n-1}},$$
which is equal to
$$\frac{p_n}{q_n}\left(1+\frac{(-1)^n}{q_n(xp_n+2q_n)+(-1)^{n-1}}\right),$$
as desired.
\end{pf}

The symmetry appearing in \propref{strange} looks very much like
2-fold symmetry, but of course is not exactly the same.  For lack of a
better name, we call it duplicating symmetry.

Given the more restrictive hypotheses of \propref{strange} (compared
to those of the Folding Lemma), one might not expect there to be any
interesting applications.  However, there are some.

Products such as $$\prod_{n=0}^\infty\left( 1 + \frac{1}{x^{\ell^n}}
\right)$$ do not have $k$-fold or duplicating symmetry.  (They do have
interesting continued fraction expansions, which can be determined
explicitly for $\ell$ even.  For the details, see \cite{products}.)
Surprisingly, the Chebyshev analogues
$$\prod_{n=0}^\infty{\left( 1 + \frac{1}{T_{\ell^n}(x)} \right)}$$
(where $T_\ell(x)$ is the $\ell$-th Chebyshev polynomial)
have iterated duplicating symmetry, provided that $\ell \equiv 2
\pmod{4}$.  We prove this in slightly greater generality.

Let $f(x)$ be a polynomial over $\Z$ of degree greater than 1 such that
$$f(x) \equiv 2x^2-1 \pmod{x(x^2-1)}$$ in $\Z[x]$.
It is easy to prove by induction that $f^i(x)^2-1$ is divisible by
$(f^0(x)+1)(f^1(x)+1)\ldots(f^i(x)+1)$.  Also, note that $f^i(x)$ divides
$f^{i+1}(x)+1$.

Now consider the product
$$P_n = \prod_{i=0}^n{\left( 1+\frac{1}{f^i(x)} \right)}.$$ The
remarks above show that its denominator is $f^n(x)$.  We now show that
it has duplicating symmetry.  We can assume (by induction) that
$P_{n-1}$ has a symmetric continued fraction, of odd length.  (Note
that $P_0 = [1,x-1,1]$.)  Now, in order to apply \propref{strange}, we
need only show that
$$\frac{f(f^{n-1}(x)) - 2f^{n-1}(x)^2+1}{f^{n-1}(x)}$$
is divisible by the numerator of $P_{n-1}$.  Note that this is divisible
by $f^{n-1}(x)^2-1$, and hence by $(f^0(x)+1)(f^1(x)+1)\ldots(f^{n-1}(x)+1)$.
This implies that it is divisible by the numerator of $P_{n-1}$, and thus
that $P_n$ has iterated duplicating symmetry.  (Note that in the
special case $f(x) = 2x^2-1$, 
we have $P_\infty = \sqrt{x^2-1}/(x-1).$  See \cite{ostrowski}.)

When $\ell \equiv 2 \pmod{4}$, the Chebyshev polynomial $T_\ell(x)$
satisfies the appropriate congruence.  Since for all $a$ and $b$,
$T_{ab}(x)=T_a(T_b(x)),$ it follows that
$$\prod_{n=0}^\infty{\left( 1 + \frac{1}{T_{\ell^n}(x)} \right)}$$
has iterated duplicating symmetry.  Also, the partial quotients have
integer coefficients, so the continued fraction is specializable.
For $\ell \equiv 0 \pmod{4}$, there appears to be a semi-specializable
continued fraction, which is somewhat, but not entirely, symmetric.
For $\ell$ odd, things are more complicated. (In these cases, the
continued fractions are almost certainly never specializable.)

\section{Open Questions and Conjectures}

Unfortunately, this paper fails to answer the following natural
questions (among others):
\begin{enumerate}
\item How do \thmref{big} and \thmref{done} generalize to products?

\item Which of the sums studied in this paper have
semi-specializable continued fractions?  What about the products?

\item There are a number of continued fraction expansions derived in
this paper that are reminiscent of 2-fold or 3-fold symmetry.  Are any
of them special cases of more general phenomena (as 2-fold symmetry is
a special case of $k$-fold symmetry)?

\item Under what conditions is the sum of two specializable continued
fractions still specializable? (At the end of \secref{appl}, a number
of cases in which this occurs are derived.  However, the derivation
doesn't seem to generalize in this direction.)

\item Does \propref{parsums} have any interesting or useful
generalizations? 
\end{enumerate}

We can, however, conjecture some partial answers to these questions:

\begin{conj}
There
exists a finite collection $\{ (a_i(x),b_i(x)) \} \subset \Z[x] \times
\Z[x]$ such that each $b_i(x)$ divides $x^2(x^2-1)^2$, with the property
that if $f(x) \in \Z[x]$ and $\deg f(x) > 1$, then all of the partial products of
$\prod_{n=0}^\infty{\left(1+1/f^n(x)\right)}$ have specializable
continued fractions if and only if for some $i$, $f(x) \equiv a_i(x)
\pmod{b_i(x)}.$ 
\end{conj}

It may be the case that $x^2(x^2-1)^2$ needs to be replaced by
another polynomial, but even in that case, the
(appropriately modified) conjecture is
presumably true.  Techniques similar to those used in the proof of
\thmref{big} will probably work.  One will presumably also arrive at a
result analogous to \thmref{done}.

It seems plausible that one can prove similar results for
semi-specializability (in the case of products or that of sums).

\begin{conj}
Let $\ell$ be a non-zero multiple of 4.  Then
$$\prod_{n=0}^\infty\left( 1+\frac{1}{T_{\ell^n}(x)}\right)$$ has a
semi-specializable continued fraction expansion.
\end{conj}

If true, this result probably arises from congruence conditions
satisfied by the Chebyshev polynomials.  The result about $x^2+x-1$
mentioned at the end of \secref{specializability} may very well also
generalize similarly.

\section*{Acknowledgements}

I am grateful to Jeffrey Shallit for reading a draft of this paper and
providing helpful comments and suggestions.

\end{document}